\documentclass[12pt]{smfart}

\usepackage[francais]{babel}
\usepackage{smfenum}
\usepackage{smfthm}
\usepackage{amssymb}
\usepackage{amscd}
\usepackage[applemac]{inputenc}

\newcommand{\Qp}{\mathbf{Q}_p}
\newcommand{\Zp}{\mathbf{Z}_p}
\newcommand{\Fp}{\mathbf{F}_p}
\newcommand{\Cp}{\mathbf{C}_p}
\newcommand{\ZZ}{\mathbf{Z}}
\newcommand{\QQ}{\mathbf{Q}}
\newcommand{\OO}{\mathcal{O}}

\newcommand{\Qph}{\mathbf{Q}_{p^h}}
\newcommand{\val}{\operatorname{val}}
\renewcommand{\phi}{\varphi}
\renewcommand{\geq}{\geqslant}
\renewcommand{\leq}{\leqslant} 
\renewcommand{\tilde}{\widetilde}
\newcommand{\vp}{\val_p}
\newcommand{\ve}{\val_{\mathrm{E}}}
\newcommand{\eps}{\varepsilon}
 
\newcommand{\bcris}{\mathbf{B}_{\mathrm{cris}}} 
\newcommand{\be}{\mathbf{B}_{\mathrm{e}}} 
\newcommand{\bmax}{\mathbf{B}_{\mathrm{max}}} 
 
\newcommand{\btp}{\widetilde{\mathbf{B}}^+_{\mathrm{rig}}} 
\newcommand{\btrig}[1]
{\widetilde{\mathbf{B}}^{\dagger #1}_{\mathrm{rig}}} 
\newcommand{\brig}[2]{\mathbf{B}^{\dagger #1}_{\mathrm{rig} #2}}

\newcommand{\bdr}{\mathbf{B}_{\mathrm{dR}}}

\newcommand{\bdag}[1]{\mathbf{B}^{\dagger #1}}

\newcommand{\btdag}[1]{\widetilde{\mathbf{B}}^{\dagger #1}}
\newcommand{\atfont}{\widetilde{\mathbf{A}}}
\newcommand{\atplus}{\widetilde{\mathbf{A}}^+}
\newcommand{\btplus}{\widetilde{\mathbf{B}}^+}

\newcommand{\bt}{\widetilde{\mathbf{B}}}
\newcommand{\et}{\widetilde{\mathbf{E}}}
\newcommand{\etplus}{\widetilde{\mathbf{E}}^+}

\newcommand{\dfont}{\mathrm{D}}

\newcommand{\dtilde}{\widetilde{\mathrm{D}}}
\newcommand{\mtilde}{\widetilde{\mathrm{M}}}

\newcommand{\fil}{\mathrm{Fil}}

\newcommand{\rg}{\mathrm{rg}}
\newcommand{\rep}{\mathrm{Rep}}
\newcommand{\bangk}{\mathcal{B}(G_K)}
\newcommand{\prcp}{\mathcal{C}(G_K)}
\newcommand{\Ob}{\mathrm{Ob}}

\newcommand{\dr}{_{\mathrm{dR}}}
\newcommand{\haut}{\mathrm{ht}}

\author{Laurent Berger}
\address{Universit\'e de Lyon \\
UMPA ENS Lyon \\
46 all\'ee d'Italie \\
69007 Lyon}
\email{laurent.berger@umpa.ens-lyon.fr}
\urladdr{www.umpa.ens-lyon.fr/\~{}lberger/}

\date{Février 2008}

\title{Presque $\mathbf{C}_p$-repr\'esentations et $(\varphi,\Gamma)$-modules}

\alttitle{Almost $\mathbf{C}_p$-representations and $(\varphi,\Gamma)$-modules}

\subjclass{11F80; 11F85; 11S15; 11S20; 11S25; 14F30}

\keywords{Théorie de Hodge $p$-adique; presque $\mathbf{C}_p$-représentations; $B$-paires; $(\varphi,\Gamma)$-modules; $(\varphi,G_K)$-modules; pentes de Frobenius; $\mathbf{B}_e$-représentations}

\numberwithin{equation}{section}

\begin{document}

\begin{abstract}
On associe deux presque $\mathbf{C}_p$-représentations à un $(\varphi,\Gamma)$-module, et on en calcule les dimensions et les hauteurs. Comme corollaire, on obtient un résultat de pleine fidélité pour les $\mathbf{B}_{\mathrm{e}}$-représentations.
\end{abstract}

\begin{altabstract}
We associate two almost $\mathbf{C}_p$-representations to a $(\varphi,\Gamma)$-module, and we compute their dimensions and heights. As a corollary, we get a full faithfulness result for $\mathbf{B}_{\mathrm{e}}$-representations.
\end{altabstract} 

\maketitle

\setcounter{tocdepth}{2}
\tableofcontents

\setlength{\baselineskip}{18pt}

\section*{Introduction}\label{secintro}

Les $(\phi,\Gamma)$-modules sont des modules sur des anneaux de séries formelles munis de structures supplémentaires qui ont été introduits dans \cite{F90} par Fontaine dans le but de décrire explicitement les représentations $p$-adiques du groupe de Galois d'un corps $p$-adique. La théorie s'est depuis pas mal développée et les $(\phi,\Gamma)$-modules que l'on considère actuellement sont des modules sur l'anneau de Robba. Grâce aux travaux de Kedlaya (voir par exemple \cite{KK04}), on a une bonne notion de pentes pour ces objets et en combinant des résultats de Fontaine (\cite{F90}), Wintenberger (\cite{FW79} et \cite{WI83}), Cherbonnier et Colmez (\cite{CC98}) et Kedlaya (\cite{KK05}) on trouve que les $(\phi,\Gamma)$-modules de pente $0$ forment une catégorie qui est naturellement équivalente à celle de toutes les représentations $p$-adiques. 

Les $B$-paires sont des objets introduits dans \cite{Ber8} qui généralisent les représentations $p$-adiques, et la catégorie des $B$-paires est alors équivalente à celle de tous les $(\phi,\Gamma)$-modules sur l'anneau de Robba. Le résultat principal de cet article est la construction de deux presque $\Cp$-représentations $X^0(W)$ et $X^1(W)$ associées à une $B$-paire $W$ ainsi que le calcul de leurs dimensions et de leurs hauteurs. Nous décrivons maintenant plus précisément les résultats de cet article. Des rappels plus fournis sont donnés aux \S\S \ref{pgbp}, \ref{prcp}.

Soit $K$ une extension finie de $\Qp$ dont on note $G_K$ le groupe de Galois absolu et $\bdr$, $\bdr^+$ et $\be=\bcris^{\phi=1}$ les anneaux construits par Fontaine. Rappelons qu'une $B$-paire est la donnée d'une $\be$-représentation $W_e$ de $G_K$ et d'une $\bdr^+$-représentation $W\dr^+$ de $G_K$ telles que $\bdr \otimes_{\bdr^+} W\dr^+ = \bdr \otimes_{\be} W_e$ (on note $W\dr$ ce dernier espace). On sait associer à toute $B$-paire un $(\phi,\Gamma)$-module $\dfont(W)$ sur l'anneau de Robba et par le théorème A de \cite{Ber8} le foncteur qui en résulte est une équivalence de catégories. Les $\phi$-modules sur l'anneau de Robba ont été étudiés par Kedlaya (voir \cite{KK05}) et on sait leur associer des pentes, et en particulier un degré $\deg(\dfont)$. Par ailleurs, le plus grand sous-module de $\dfont$ de pentes $\leq 0$ est bien défini, et on le note $\dfont_{\leq 0}$; on note aussi $\dfont_{>0} = \dfont / \dfont_{\leq 0}$. 

Les presque $\Cp$-représentations sont des objets définis et étudiés par Fontaine dans \cite{F03}. Ce sont des espaces de Banach $X$ munis d'une action linéaire et continue de $G_K$ tels qu'il existe $d\geq 0$ et des représentations $p$-adiques de dimensions finies $V_1$ et $V_2$ telles que $X/V_1 = \Cp^d / V_2$. Les presque $\Cp$-représentations ont une dimension $\dim_{\prcp}(X)$ (l'entier $d$ ci-dessus) et une hauteur $\haut(X)$ (qui vaut $\dim(V_1)-\dim(V_2)$). 

Rappelons enfin qu'en plus de l'inclusion $\be \subset \bdr$, ces deux anneaux sont reliés par la suite exacte fondamentale $0 \to \Qp \to \be \to \bdr/\bdr^+ \to 0$.
Si $W$ est une $B$-paire, alors on pose $X^0(W) = \ker(W_e \to W\dr / W\dr^+)$ et $X^1(W) = \mathrm{coker}(W_e \to W\dr / W\dr^+)$. Le résultat principal de notre article est alors le suivant (voir le théorème \ref{xipcp} qui donne un résultat plus précis) :
\begin{enonce*}{Théorème A}
Si $W$ est une $B$-paire, alors  $X^0(W)$ et $X^1(W)$ sont deux presque $\Cp$-représentations. Si de plus $\dfont=\dfont(W)$, alors :
\begin{enumerate}
\item $\dim_{\prcp}X^0(W) = -\deg(\dfont_{\leq 0})$ et $\haut(X^0(W)) = \rg(\dfont_{\leq 0})$;
\item $\dim_{\prcp}X^1(W) = \deg(\dfont_{> 0})$ et $\haut(X^1(W)) = -\rg(\dfont_{> 0})$. 
\end{enumerate}
\end{enonce*}

Une $\be$-représentation est de manière naturelle une limite inductive d'espaces de Banach munis d'une action linéaire et continue de $G_K$. Par ailleurs, toute $\be$-représentation peut être vue comme le $W_e$ d'une $B$-paire et en utilisant le théorème A, on trouve le résultat suivant (c'est le théorème \ref{ffbe}), qui répond par l'affirmative à une question de Fontaine (voir la remarque en bas de la page 375 de \cite{F03}).

\begin{enonce*}{Théorème B}
Le foncteur d'oubli de la catégorie des $\be$-représentations de $G_K$ vers la catégorie des $\Qp$-espaces vectoriels topologiques avec action linéaire et continue de $G_K$ est pleinement fidèle.
\end{enonce*}

Le théorème A est dans le prolongement naturel des calculs de \cite{Ber8}. Notons toutefois que le fait que $X^0(W)$ et $X^1(W)$ sont des presque $\Cp$-représentations ainsi que le théorème B s'inscrivent naturellement dans un projet ambitieux ayant pour but de faire le lien entre $B$-paires, presque $\Cp$-représentations et espaces de Banach-Colmez (Fontaine et Plût, travail en cours). Je remercie en particulier Jean-Marc Fontaine pour des discussions éclairantes en rapport avec le contenu de cet article.

La plan de l'article est assez naturel. Dans le \S\ref{pgbp}, on fait quelques rappels et compléments sur les $B$-paires, les $(\phi,G_K)$-modules et les $(\phi,\Gamma)$-modules et dans le \S\ref{prcp} on fait de même pour les presque $\Cp$-représentations. Le \S\ref{subalcp} est consacré à la démonstration du théorème A et le \S\ref{parffbe} à celle du théorème $B$. Enfin, dans un court appendice on donne la démonstration d'un énoncé utilisé implicitement dans \cite{Ber5} et \cite{Ber8}.
 
\section{$B$-paires, $(\phi,G_K)$-modules et $(\phi,\Gamma)$-modules}\label{pgbp}

Commen\c{c}ons par faire des rappels tr\`es 
succints sur les d\'efinitions
(donn\'ees dans \cite{F3} par exemple)
des divers anneaux que nous utilisons dans cet article. Rappelons
que $\etplus = \varprojlim_{x \mapsto x^p} \OO_{\Cp}$ est
un anneau de caract\'eristique $p$, complet pour la valuation
$\ve$ d\'efinie par $\ve(x)= 
\vp(x^{(0)})$ et qui contient un \'el\'ement $\eps$ tel que
$\eps^{(n)}$ est une racine primitive $p^n$-i\`eme de l'unit\'e.
On fixe un tel $\eps$ dans tout l'article.
L'anneau $\et = \etplus[1/(\eps-1)]$ est alors un corps qui 
contient comme sous-corps dense la cl\^oture alg\'ebrique
de $\Fp(\!(\eps-1)\!)$. On pose $\atplus=W(\etplus)$ et 
$\atfont=W(\et)$ ainsi que $\btplus=\atplus[1/p]$ et 
$\bt=\atfont[1/p]$. 
L'application $\theta : \btplus \to \Cp$ qui \`a
$\sum_{k \gg -\infty} p^k [x_k]$ associe 
$\sum_{k \gg -\infty} p^k x_k^{(0)}$ est un morphisme
d'anneaux surjectif et $\bdr^+$ est le compl\'et\'e de $\btplus$
pour la topologie $\ker(\theta)$-adique, ce qui en fait un espace
topologique de Fr\'echet. On pose $X=[\eps]-1 \in \atplus$ et
$t=\log(1+X) \in \bdr^+$ et on d\'efinit $\bdr$ par $\bdr=\bdr^+[1/t]$.
Soit $\tilde{p} \in \etplus$ un \'el\'ement tel que $p^{(0)}=p$.
L'anneau $\bmax^+$ est l'ensemble des s\'eries $\sum_{n \geq 0} b_n 
([\tilde{p}]/p)^n$ o\`u $b_n \in \btplus$ et $b_n \to 0$ quand $n \to \infty$
ce qui en fait un sous-anneau de $\bdr^+$ 
muni en plus d'un Frobenius $\phi$ qui est injectif, mais pas surjectif.
On pose $\btp  = \cap_{n\geq 0} \phi^n(\bmax^+)$ ce qui
en fait un sous-anneau de $\bmax^+$ sur lequel $\phi$ est bijectif.
Remarquons que l'on travaille souvent avec $\bcris$ plut\^ot que
$\bmax$ mais le fait de pr\'ef\'erer $\bmax$ ne change rien aux
r\'esultats et est plus agr\'eable pour des raisons techniques.

Rappelons que les anneaux $\bmax$
et $\bdr$ sont reli\'es, en plus de l'inclusion $\bmax \subset \bdr$, par
la suite exacte fondamentale $0 \to \Qp \to \bmax^{\phi=1} 
\to \bdr/\bdr^+ \to 0$.
Ce sont ces anneaux que l'on utilise en th\'eorie de Hodge $p$-adique. 
Le point de d\'epart de la th\'eorie des $(\phi,\Gamma)$-modules 
sur l'anneau de Robba est la construction d'anneaux
interm\'ediaires entre $\btplus$ et $\bt$. Si $r>0$, soit $\btdag{,r}$
l'ensemble des $x = \sum_{k \gg -\infty} p^k [x_k] \in \bt$ tels que
$\ve(x_k)+k \cdot pr/(p-1)$ tend vers $+\infty$ quand $k$ augmente. 
On pose $\btdag{} = \cup_{r > 0} \btdag{,r}$, c'est le corps
des \'el\'ements surconvergents, d\'efini dans \cite{CC98}.
L'anneau $\btrig{} =  \cup_{r > 0} \btrig{,r}$ 
d\'efini dans \cite[\S 2.3]{Ber1} est en quelque sorte
la somme de $\btp$ et $\btdag{}$; de fait, on a une suite exacte (d'anneaux et d'espaces de Fr\'echet) $0 \to \btplus \to \btp \oplus \btdag{,r} \to \btrig{,r} \to 0$.

Rappelons que $K_0=W(k)[1/p]$; pour $1 \leq n \leq +\infty$, 
on pose $K_n=K(\mu_{p^n})$ et $H_K = \mathrm{Gal}
(\overline{K}/K_\infty)$ et $\Gamma_K=G_K/H_K$.
Si $R$ est un anneau muni d'une action de $G_K$ (c'est le cas
pour tous ceux que nous consid\'erons), on note $R_K = R^{H_K}$.
L'anneau $\btrig{,r}$ contient l'ensemble des s\'eries 
$f(X) = \sum_{k \in \ZZ} f_k X^k$ avec $f_k \in K_0$ telles que
$f(X)$ converge sur $\{ p^{-1/r} \leq |X| < 1 \}$. Cet anneau est
not\'e $\brig{,r}{,K_0}$. Si $K$ est une extension finie de $K_0$, il
lui correspond par la th\'eorie du corps de normes 
(cf.\ \cite{FW79} et \cite{WI83})
une extension finie
$\brig{,r}{,K}$ qui s'identifie (si $r$ est assez grand) 
\`a l'ensemble des s\'eries 
$f(X_K) = \sum_{k \in \ZZ} f_k X_K^k$ avec $f_k \in K'_0$ telles que
$f(X_K)$ converge sur $\{ p^{-1/e r} \leq |X_K| < 1 \}$ o\`u 
$X_K$ est un certain \'el\'ement de $\btdag{}_K$ et $K'_0$ est
la plus grande sous-extension non ramifi\'ee de $K_0$ dans 
$K_\infty$ et $e=[K_\infty : K_0(\mu_{p^\infty})]$.
On pose $\brig{}{,K} = \cup_{r > 0} \brig{,r}{,K}$ et
$\bdag{,r}_K = \brig{,r}{,K} \cap \btdag{}$ et  $\bdag{}_K =
\cup_{r > 0}  \bdag{,r}_K$.
Les anneaux $\btrig{}$ et $\brig{}{,K}$ 
co\"{\i}ncident avec 
les anneaux $\Gamma^{\mathrm{alg}}_{\mathrm{an,con}}$
et $\Gamma_{\mathrm{an,con}}$ d\'efinis
dans \cite[\S 2.2]{KK05} (cf.\ en particulier la convention 2.2.16 et la
remarque 2.4.13 de \cite{KK05}). 

Rappelons que l'on pose $Q_1=((1+X)^p-1)/X$ et $Q_n=\phi^{n-1}(Q_1)$ pour $n \geq 1$ et que si $r>0$, alors $n(r)$ est le plus petit entier $n$ tel que $p^{n-1}(p-1) \geq r$. Enfin, l'injection $\phi^{-n} : \btplus \to \bdr^+$ se prolonge par continuité à $\iota_n : \btrig{,p^{n-1}(p-1)} \to \bdr^+$. 

\begin{lemm}\label{divpart}
Si $x \in \btrig{,r}$ a la propriété que pour tout $n \geq n(r)$, on a $\theta \circ \iota_n(x) = 0$, alors $x \in t \cdot \btrig{,r}$.
\end{lemm}

\begin{proof}
Soit $J$ l'idéal de $\btrig{,r}$ constitué des $x$ vérifiant l'hypothèse du lemme. Comme les applications $\theta \circ \iota_n : \btrig{,r} \to \Cp$ sont continues, l'idéal $J$ est fermé et par le théorème 2.9.6 de \cite{KK05}, il est donc principal. On a manifestement $t \in J$, et l'idéal $J$ est donc engendré par un élément $f$ qui divise $t$. Par la proposition 2.17 de \cite{Ber1}, le noyau de $\theta \circ \iota_n : \btrig{,r} \to \Cp$ est engendré par $Q_n$ et on en déduit que $Q_n$ divise $f$ quel que soit $n \geq n(r)$ et comme dans $\btrig{,r}$ les idéaux engendrés par $t$ et par $\prod_{n \geq n(r)} (Q_n / p)$ co\"{\i}ncident, on a bien $J = t \cdot \btrig{,r}$.
\end{proof}

Rappelons qu'un $\phi$-module sur $\brig{}{,K}$ est un 
$\brig{}{,K}$-module libre $\dfont$ muni d'un Frobenius $\phi$
tel que $\phi^*(\dfont)=\dfont$. Un $(\phi,\Gamma)$-module
est un $\phi$-module muni d'une action semi-lin\'eaire 
et continue de $\Gamma_K$ qui commute \`a $\phi$.
Les $B$-paires sont des objets qui ont été définis dans \cite{Ber8}. Rappelons qu'une $B$-paire $W = (W_e,W\dr^+)$ est la donnée d'une $\be$-représentation $W_e$ et d'une $\bdr^+$-représentation $W\dr^+$ qui satisfont $\bdr \otimes_{\be} W_e = \bdr \otimes_{\bdr^+} W\dr^+$. Dans \cite{Ber8}, nous avons construit une équivalence de catégories entre la catégorie des $B$-paires et celle des $(\phi,\Gamma)$-modules sur $\brig{}{,K}$. Nous allons la préciser ci-dessous.

Appelons $(\phi,G_K)$-module sur $\btrig{}$ la donnée d'un $\phi$-module $\dtilde$ sur $\btrig{}$ muni d'une action de $G_K$ qui est semi-linéaire et continue et qui commute à $\phi$. On a un foncteur évident $\dfont \mapsto \dtilde = \btrig{} \otimes_{\brig{}{,K}} \dfont$ de la catégorie des $(\phi,\Gamma)$-modules sur $\brig{}{,K}$ vers celle des $(\phi,G_K)$-modules sur $\btrig{}$ et on peut par ailleurs associer à tout $(\phi,G_K)$-module $\dtilde$ sur $\btrig{}$ une $B$-paire $W(\dtilde)$ en copiant la recette de la proposition 2.2.6 de \cite{Ber8}. Le théorème 2.2.7 de \cite{Ber8} s'étend alors immédiatement en le résultat suivant :

\begin{theo}\label{plus}
Les foncteurs $\dfont \mapsto \dtilde$ et $\dtilde \mapsto W(\dtilde)$ réalisent des équivalences de catégories entre les trois catégories suivantes :
\begin{enumerate}
\item les $(\phi,\Gamma)$-modules sur $\brig{}{,K}$;
\item les $(\phi,G_K)$-modules sur $\btrig{}$;
\item les $B$-paires.
\end{enumerate}
\end{theo}

Remarquons que si $\dfont(W)$ est le $(\phi,\Gamma)$-module associé à une $B$-paire $W= (W_e,W\dr^+)$, alors le $(\phi,\Gamma)$-module associé à $W= (W_e, t^j  W\dr^+)$ est $t^j  \dfont(W)$.

Kedlaya a étudié dans \cite{KK05} les $\phi$-modules sur $\btrig{}$ et il a notamment dégagé la notion de pentes pour ces objets (il a en fait un cadre commun pour l'étude des $\phi$-modules sur $\brig{}{,K}$ et sur $\btrig{}$; pour des raisons pratiques, nous travaillons avec des modules sur $\btrig{}$ dans cet article). 

Si $h \geq 1$ et $a \in \ZZ$, disons qu'un $\phi$-module $\dtilde$ sur $\btrig{}$ est isocline de pente $a/h$ si c'est une somme directe des modules élémentaires $M_{a,h}$ définis dans le \S 4.1 de \cite{KK05}. On peut appliquer le même raisonnement que dans la démonstration du théorème 3.2.3 de \cite{Ber8} pour en déduire que tout $\phi$-module $\dtilde$ sur $\btrig{}$ qui est isocline de pente $a/h$ s'écrit $\dtilde = \btrig{} \otimes_{\Qph} V_{a,h}$ où $V_{a,h} = \dtilde^{\phi^h = p^a}$ est un $\Qph$-espace vectoriel de dimension finie muni d'un frobenius semi-linéaire vérifiant $\phi^h=p^a$. Si de plus $\dtilde$ est un $(\phi,G_K)$-module sur $\btrig{}$, alors $V_{a,h}$ est stable par $G_K$ et donc $V_{a,h}$ appartient à la catégorie $\rep(a,h)$ définie dans \cite[\S 3.2]{Ber8} et on a :

\begin{prop}\label{modisoc}
Le foncteur $V_{a,h} \mapsto \btrig{} \otimes_{\Qph} V_{a,h}$ réalise une équivalence de catégories entre la catégorie $\rep(a,h)$ et celle des $(\phi,G_K)$-modules sur $\btrig{}$ qui sont isoclines de pente $a/h$.
\end{prop}

En ce qui concerne les $(\phi,G_K)$-modules non nécessairement isoclines, Kedlaya a montré le résultat ci-dessous (voir \cite[\S 4.5]{KK05}).

\begin{theo}\label{dmked}
Si $\dtilde$ est un $\phi$-module sur $\btrig{}$, alors il existe une unique filtration $0 = \dtilde_0 \subset \dtilde_1 \subset \cdots \subset \dtilde_\ell = \dtilde$ par des sous-$\phi$-modules saturés, telle que :
\begin{enumerate}
\item pour tout $1 \leq i \leq \ell$, le quotient $\dtilde_i / \dtilde_{i-1}$ est isocline;
\item si l'on appelle $s_i$ la pente de $\dtilde_i / \dtilde_{i-1}$, alors $s_1 < s_2 < \cdots < s_\ell$.
\end{enumerate}
De plus, cette filtration est scindée (décomposition de Dieudonné-Manin), mais non canoniquement.
\end{theo}
 
Le fait que cette filtration est canonique implique que si $\dtilde$ est un $(\phi,G_K)$-module sur $\btrig{}$, alors chacun des crans de la filtration est stable par $G_K$ et est donc aussi un $(\phi,G_K)$-module. En revanche, il n'est en général pas possible de scinder la filtration de manière compatible à l'action de $G_K$ (ceci dit, c'est sans doute possible si $\dtilde$ est de de Rham, voir la décomposition de Dieudonné-Manin de \cite{PC03}).

Les rationnels $s_1,\cdots,s_\ell$ construits dans le théorème \ref{dmked} sont par définition les pentes de $\dtilde$. Si $s \in \QQ$, alors $\dtilde_{\leq s}$ et $\dtilde_{< s}$, les sous-modules maximaux de $\dtilde$ de pentes $\leq s$ et $<s$ sont bien définis et ce sont des $(\phi,G_K)$-modules. On note aussi $\dtilde_{>s} = \dtilde / \dtilde_{\leq s}$ et $\dtilde_{\geq s} = \dtilde / \dtilde_{< s}$.

Pour terminer, nous aurons besoin de la proposition ci-dessous, qui généralise un peu la proposition 4.1.3 de \cite{KK05}.

\begin{prop}\label{hiphi}
Soit $\dtilde$ un $(\phi,G_K)$-module sur $\btrig{}$ et $h \geq 1$ et $a\in \ZZ$.
\begin{enumerate}
\item Si les pentes de $\dtilde$ sont $> a/h$, alors $\dtilde^{\phi^h=p^a} = 0$; 
\item si les pentes de $\dtilde$ sont $\leq a/h$, alors $\dtilde/(\phi^h - p^a) = 0$.
\end{enumerate}
\end{prop}

\begin{proof}
Cela suit du corollaire 4.1.4 de \cite{KK05} puisque l'on a d'une part $\dtilde^{\phi^h=p^a} = \mathrm{Hom}(M_{a,h},\dtilde)$ et d'autre part $\dtilde/(\phi^h - p^a) = \mathrm{Ext}^1(M_{a,h},\dtilde)$.
\end{proof}

\section{Rappels et compléments sur les presque $\Cp$-représentations}\label{prcp}
Soit $\bangk$ la catégorie dont les objets sont les $\Qp$-espaces de Banach munis d'une action linéaire et continue de $G_K$ et dont les morphismes sont les applications linéaires continues et $G_K$-equivariantes, et soit $\prcp$ la sous-catégorie de $\bangk$ constituée des presque $\Cp$-représentations de $G_K$ définies et étudiées dans \cite{F03}. Nous rappelons ici quelques uns des résultats à leur sujet qui nous serviront dans la suite de cet article. 

Par le théorème B de \cite{F03}, la catégorie $\prcp$ est une sous-catégorie stricte de $\bangk$ (c'est-à-dire une sous-catégorie strictement pleine, stable par somme directe, et telle que tout morphisme est strict et a son noyau et son conoyau dans $\prcp$), et  $\prcp$ contient toutes les représentations $p$-adiques de dimension finie et toutes les $\bdr^+$-représentations de longueur finie. 

De plus, il existe deux fonctions, la dimension $\dim_{\prcp} : \Ob\ \prcp \to \ZZ_{\geq 0}$ et la hauteur $\haut : \Ob\ \prcp \to \ZZ$, qui sont additives et caractérisées par $\dim_{\prcp}(\Cp)=1$ et $\haut(\Cp)=0$ et $\dim_{\prcp}(\Qp)=0$ et $\haut(\Qp)=1$. On a par exemple $\dim_{\prcp}(\bdr^+/t^a)=a$ et $\haut(\bdr^+/t^a)=0$. Une presque $\Cp$-représentation de dimension et de hauteur nulles est elle-même nulle.

\begin{lemm}\label{btrigprcp}
Si $h \geq 1$ et $a \in \ZZ$, alors $(\btrig{})^{\phi^h=p^a}$ est une presque $\Cp$-représen\-tation, qui est de dimension $a$ et de hauteur $h$ si $a \geq 0$ et qui est nulle si $a < 0$.
\end{lemm}

\begin{proof}
La proposition 2.10 de \cite{PC03} implique que si $a\geq 0$, alors on a une suite exacte :
\[ 0 \to \Qph \cdot t_h^a \to (\btrig{})^{\phi^h=p^a} \to \bdr^+ / t^a \to 0, \] où $t_h$ est l'élément construit dans le \S 2.4 de \cite{PC03}, et donc $(\btrig{})^{\phi^h=p^a}$ est bien une presque $\Cp$-représentation et sa dimension se lit sur la suite exacte.
\end{proof}

\begin{prop}\label{prisoc0}
Si $h,k \geq 1$ et $a,b \in \ZZ$, et si $V_{a,h} \in \rep(a,h)$ est de dimension $d$ sur $\Qph$, alors $(\btrig{} \otimes_{\Qph} V_{a,h})^{\phi^k=p^b}$ est une presque $\Cp$-représentation, nulle si $a/h > b/k$, et de dimension $d \cdot (b-ak/h)$ et de hauteur $dk$ sinon.
\end{prop}

\begin{proof}
Un calcul immédiat montre que :
\[ X = (\btrig{} \otimes_{\Qph} V_{a,h})^{\phi^k=p^b} \subset  Y = (\btrig{})^{\phi^{kh}=p^{bh-ak}} \otimes_{\Qph} V_{a,h} \] 
et le lemme \ref{btrigprcp} montre que $Y$ est une presque $\Cp$-représentation, qui est nulle si $a/h > b/k$. L'espace $Y$ est muni de l'opérateur $\phi$ qui est continu et commute à l'action de $G_K$. Comme $\prcp$ est une sous-catégorie pleine de $\bangk$, l'espace $X=Y^{\phi^k-p^b}$ est lui-même un objet de $\prcp$.

Observons que si $a/h \leq b/k$, alors $\dim_{\prcp}(Y)=d(bh-ak)$ et $\haut(Y)=dkh$ par le lemme \ref{btrigprcp} et comme on a un isomorphisme $Y \simeq \QQ_{p^{kh}} \otimes_{\QQ_{p^k}}  X$ (voir par exemple la proposition 3.18 de \cite{PC03}), on trouve que $\dim_{\prcp}(X)=d(bh-ak)/h$ et que $\haut(X)=dk$ (notons que par la remarque 3.2.2 de \cite{Ber8}, $d$ est divisible par $h$).
\end{proof}

\begin{prop}\label{prisoc1}
Si $h \geq 1$ et $a \in \ZZ$, et si $V_{a,h} \in \rep(a,h)$ est de dimension $d$ sur $\Qph$, alors :
\[ X^1(V_{a,h}) = \frac{\bdr \otimes_{\Qph} V_{a,h}} {\bdr^+ \otimes_{\Qph} V_{a,h} + (\btrig{}[1/t] \otimes_{\Qph} V_{a,h})^{\phi=1}} \] 
est une presque $\Cp$-représentation, nulle si $a \leq 0$ et de dimension $da/h$ et de hauteur $-d$ sinon.
\end{prop}

\begin{proof}
Si $a \leq 0$ et si $\ell \geq 0$, alors un petit calcul et la proposition 2.11 de \cite{PC03} impliquent qu'on a un morceau de suite exacte :
\[ 0 \to (\btrig{} \otimes_{\Qph} V_{a,h})^{\phi=1} \to (t^{-\ell} \btrig{} \otimes_{\Qph} V_{a,h})^{\phi=1} \to (t^{-\ell} \bdr^+ / \bdr^+) \otimes_{\Qph} V_{a,h}, \]
et on a $(t^{-\ell} \btrig{} \otimes_{\Qph} V_{a,h})^{\phi=1} = t^{-\ell}  (\btrig{} \otimes_{\Qph} V_{a,h})^{\phi=p^{\ell}}$. Par la proposition \ref{prisoc0}, les deux premiers termes sont de dimensions $-da/h$ et $d \cdot (\ell-a/h)$ et de hauteurs $d$ et $d$ ce qui fait que la flèche de droite est en fait surjective et donc :
\[ \bdr^+ \otimes_{\Qph} V_{a,h} + (t^{-\ell} \btrig{} \otimes_{\Qph} V_{a,h})^{\phi=1} = t^{-\ell} \bdr^+ \otimes_{\Qph} V_{a,h}. \]

On déduit des calculs précédents que si $j \geq 0$ vérifie $j - a/h \geq 0$, alors pour tout $\ell \gg 0$,  on a :
\[ t^{-j} \bdr^+ \otimes_{\Qph} V_{a,h} + (t^{-\ell} \btrig{} \otimes_{\Qph} V_{a,h})^{\phi=1} = t^{-\ell} \bdr^+ \otimes_{\Qph} V_{a,h}, \]
et donc :
\[ t^{-j} \bdr^+ \otimes_{\Qph} V_{a,h} + (\btrig{}[1/t] \otimes_{\Qph} V_{a,h})^{\phi=1} = \bdr \otimes_{\Qph} V_{a,h}, \]
ce qui fait que l'application :
\[ \frac{t^{-j} \bdr^+ \otimes_{\Qph} V_{a,h}} {\bdr^+ \otimes_{\Qph} V_{a,h} + (t^{-j} \btrig{} \otimes_{\Qph} V_{a,h})^{\phi=1}} \to \frac{\bdr \otimes_{\Qph} V_{a,h}} {\bdr^+ \otimes_{\Qph} V_{a,h} + (\btrig{}[1/t] \otimes_{\Qph} V_{a,h})^{\phi=1}} \]
est un isomorphisme et que $X^1(V_{a,h})$ s'identifie au conoyau de l'application :
\[ (t^{-j} \btrig{} \otimes_{\Qph} V_{a,h})^{\phi = 1} \to (t^{-j} \bdr^+/ \bdr^+) \otimes_{\Qph} V_{a,h}. \]
Comme il s'agit d'un morphisme de presque $\Cp$-représentations, ce conoyau est aussi une presque $\Cp$-représentation, et si $a \leq 0$, alors on peut prendre $j=0$ et $X^1(V_{a,h})$ est nul. Si $a \geq 1$, alors $(\btrig{} \otimes_{\Qph} V_{a,h})^{\phi=1} = 0$ et on a la suite exacte :
\[ 0 \to (t^{-j} \btrig{} \otimes_{\Qph} V_{a,h})^{\phi=1} \\
\to (t^{-j} \bdr^+ / \bdr^+) \otimes_{\Qph} V_{a,h} \to X^1(V_{a,h}) \to 0, \]
ce qui fait que $\dim_{\prcp}X^1(V_{a,h}) = dj - d(j-a/h) = d a/h$ et que $\haut(X^1(V_{a,h})) = -d$.  
\end{proof}

\begin{prop}\label{extpcp}
Si $V$ est une représentation $p$-adique de $G_K$, si $S$ est une presque $\Cp$-représentation de $G_K$ et si $E \in \bangk$ est une extension de $V$ par $S$, alors :
\begin{enumerate}
\item $E$ est elle-même une presque $\Cp$-représentation;
\item il existe un $\Qp$-espace vectoriel de dimension finie $X \subset E$ stable par $G_K$ qui s'envoie surjectivement sur $V$.
\end{enumerate}
\end{prop}

\begin{proof}
Le (1) suit de la proposition 6.7 de \cite{F03}. Par le corollaire du \S 5.2 de \cite{F03}, la suite exacte $0 \to S \to E \to V \to 0$ est presque scindée ce qui veut dire qu'il existe une représentation de dimension finie $W \subset S$ telle que la suite $0 \to S/W \to E/W \to V \to 0$ est scindée, disons par $s : V \to E/W$. Pour montrer le (2), il suffit de prendre pour $X$ l'image inverse de $s(V)$ dans $E$. 
\end{proof}

\section{Les presque $\Cp$-représentations associées aux $B$-paires}\label{subalcp}

Si $W$ est une $B$-paire, alors on pose $X^0(W) =  W_e \cap
W\dr^+ \subset W\dr$ et $X^1(W) = W\dr / ( W_e + 
W\dr^+)$. Ces espaces sont donc le noyau et le conoyau de l'application naturelle $W_e \to W\dr/W\dr^+$. L'objet de ce paragraphe est de montrer le théorème ci-dessous.

\begin{theo}\label{xipcp}
Si $W$ est une $B$-paire, et si $\dtilde=\dtilde(W)$, alors :
\begin{enumerate}
\item $X^0(W) \simeq \dtilde^{\phi=1}$ et $X^1(W) \simeq \dtilde/(1-\phi)$;
\item $X^0(W)$ et $X^1(W)$ sont deux presque $\Cp$-représentations;
\item $\dim_{\prcp}X^0(W) = -\deg(\dtilde_{\leq 0})$ et $\haut(X^0(W)) = \rg(\dtilde_{\leq 0})$;
\item $\dim_{\prcp}X^1(W) = \deg(\dtilde_{> 0})$ et $\haut(X^1(W)) = -\rg(\dtilde_{> 0})$. 
\end{enumerate}
\end{theo}

Avant de pouvoir montrer ce théorème, nous allons devoir établir quelques résultats intermédiaires. Dans la suite, il est plus commode de travailler avec des $(\phi,G_K)$-modules sur $\btrig{}$ qu'avec des $B$-paires et on pose donc $X^i(\dtilde)=X^i(W(\dtilde))$ pour $i=0$, $1$.

\begin{prop}\label{xzdtphi}
Si $\dtilde$ est un $(\phi,G_K)$-module sur $\btrig{}$, alors $X^0(\dtilde) = \dtilde^{\phi=1}$.
\end{prop}

\begin{proof}
C'est l'occasion de rappeler que par la définition donnée dans la proposition 2.2.6 de \cite{Ber8}, on a $W_e(\dtilde) = (\btrig{}[1/t] \otimes_{\btrig{}} \dtilde)^{\phi=1}$ et $W\dr^+(\dtilde) = \bdr^+ \otimes_{\btrig{r_n}}^{\iota_n} \dtilde^{r_n}$ pour $n \gg 0$. Si $\{f_i\}_{i=1}^d$ est une base de $\dtilde$ et si $y = \sum_{i=1}^d b_i \otimes f_i \in W_e$ a la propriété que :
\[ y = \phi^{-n}(y) = \sum_{i=1}^d \phi^{-n}(b_i) \otimes \phi^{-n} (f_i) \in \bdr^+ \otimes_{\btrig{r_n}}^{\iota_n} \dtilde^{r_n} \] pour tout $n \gg 0$, c'est donc que $\phi^{-n}(b_i) \in \bdr^+$ pour tout $n \gg 0$ et on est donc ramené à montrer que si $b \in \btrig{}[1/t]$ vérifie $\phi^{-n}(b) \in \bdr^+$ pour tout $n \gg 0$, alors $b \in \btrig{}$. Si $h \geq 1$ est tel que $t^h b \in \btrig{}$, alors $\theta \circ \iota_n (t^h b)= 0$ pour tout $n \gg 0$ et par le lemme \ref{divpart}, $t^h b \in t \cdot \btrig{}$ ce qui fait que l'on peut diminuer $h$ de $1$ et finalement que $b \in \btrig{}$.
\end{proof}

\begin{prop}\label{xundtsurphi}
Si $\dtilde$ est un $(\phi,G_K)$-module sur $\btrig{}$, alors $X^1(\dtilde) = \dtilde/(1-\phi)$.
\end{prop}

\begin{proof}
Si l'on a une suite exacte du type $0 \to \dtilde' \to \dtilde \to \dtilde'' \to 0$, alors le lemme du serpent appliqué au diagramme :

{\smaller \[ \begin{CD} 
0 @>>> W_e(\dtilde') \oplus W\dr^+(\dtilde') @>>> W_e(\dtilde) \oplus W\dr^+(\dtilde) @>>> W_e(\dtilde'') \oplus W\dr^+(\dtilde'') @>>> 0 \\
 && @VVV @VVV @VVV \\
 0 @>>> W\dr(\dtilde') @>>> W\dr(\dtilde) @>>> W\dr(\dtilde'') @>>> 0 
\end{CD} \]}

nous donne une suite exacte :
\begin{equation} \label{longxi}
0 \to X^0(\dtilde') \to X^0(\dtilde) \to X^0(\dtilde'') \to X^1(\dtilde') \to X^1(\dtilde) \to X^1(\dtilde'') \to 0. 
\end{equation}
Si $\dtilde$ est isocline de pente $\leq 0$, alors les propositions \ref{modisoc} et \ref{prisoc1} montrent que $X^1(\dtilde)=0$. On déduit du théorème \ref{dmked} et de la suite exacte \ref{longxi} que si $\dtilde$ est un $(\phi,G_K)$-module, alors $X^1(\dtilde_{\leq 0}) = 0$ et donc que l'application $X^1(\dtilde) \to X^1(\dtilde_{>0})$ est un ismorphisme.

Le lemme du serpent appliqué au diagramme :
\[ \begin{CD}
0 @>>> \dtilde' @>>> \dtilde @>>> \dtilde'' @>>> 0 \\
 && @VV{1-\phi}V @VV{1-\phi}V @VV{1-\phi}V \\
0 @>>> \dtilde' @>>> \dtilde @>>> \dtilde'' @>>> 0
\end{CD} \]
nous donne une suite exacte :
\begin{equation} \label{longdphi}
0 \to (\dtilde')^{\phi=1} \to \dtilde^{\phi=1} \to (\dtilde'')^{\phi=1} \to
\dtilde'/(1-\phi) \to \dtilde/(1-\phi) \to \dtilde''/(1-\phi) \to 0. 
\end{equation}
On en déduit comme ci-dessus, en utilisant cette fois le (2) de la proposition \ref{hiphi}, que si $\dtilde$ est un $(\phi,G_K)$-module, alors l'application $\dtilde/(1-\phi) \to\dtilde_{>0}/(1-\phi)$ est un isomorphisme.

Il suffit donc de montrer la proposition sous l'hypothèse supplémentaire que les pentes de $\dtilde$ sont $>0$. La proposition A.4 de \cite{KK08} nous dit qu'il existe alors deux $(\phi,G_K)$-modules $\dtilde_{\mathrm{et}}$ et $\mtilde$ tels que $\dtilde_{\mathrm{et}}$ est étale et tels que l'on a une suite exacte :
\[ 0 \to \dtilde \to \dtilde_{\mathrm{et}} \to \mtilde \to 0. \] 
En utilisant la proposition \ref{xzdtphi} et les deux suites exactes ci-dessus, ainsi que le fait que $X^1(\dtilde_{\mathrm{et}})=0$ et $\dtilde_{\mathrm{et}}/(1-\phi)=0$, on trouve que l'on a :
\[ \begin{CD} 0 @>>> X^0(\dtilde) @>>> X^0(\dtilde_{\mathrm{et}}) @>>> X^0(\mtilde) @>>> X^1(\dtilde) @>>>  0 \\
&& @| @| @| \\  
0 @>>> \dtilde^{\phi=1} @>>> (\dtilde_{\mathrm{et}})^{\phi=1} @>>> (\mtilde)^{\phi=1} @>>>
\dtilde/(1-\phi) @>>> 0, 
\end{CD} \]
ce qui fait que l'on a bien $X^1(\dtilde) \simeq \dtilde/(1-\phi)$.
\end{proof}

\begin{rema}\label{dsurphi}
L'isomorphisme $\alpha : \dtilde/ (1-\phi) \to X^1(\dtilde)$ peut être décrit explicitement de la manière suivante. Par le théorème \ref{dmked} ci-dessus et le corollaire 1.1.6 de \cite{Ber8}, l'application $1-\phi : \dtilde[1/t] \to \dtilde[1/t]$ est surjective et si $x \in \dtilde$, il existe donc $y \in \dtilde[1/t]$ tel que $(1-\phi)y=x$. Si $n \gg 0$, alors on a $\phi^{-n}(y) \in W\dr(\dtilde)$ et $\alpha(x)$ est par définition l'image de $\phi^{-n}(y)$ dans $X^1(\dtilde)$. 

Vérifions que $\alpha$ est bien définie. Si l'on avait choisi $y' \in \dtilde[1/t]$ tel que $(1-\phi)y' = x$, alors $y-y' \in (\dtilde[1/t])^{\phi=1}$ et alors $\phi^{-n}(y)$ et $\phi^{-n}(y')$ ont même image modulo $W_e(\dtilde)$ et donc $\alpha(x)$ ne dépend pas du choix de $y$. Ensuite, on a $\phi^{-(n+1)}(y) - \phi^{-n}(y) = \phi^{-(n+1)}(1-\phi)y = \phi^{-(n+1)}(x) \in W\dr^+(\dtilde)$ et donc $\alpha(x)$ ne dépend pas du choix de $n \gg 0$ ce qui fait que $\alpha(x)$ est bien défini.

On peut voir directement que $\alpha$ est injective. Si $\alpha(x)=0$, c'est que pour tout $n \gg 0$ on peut écrire $\phi^{-n}(y) = z_n + w_n$ avec $z_n \in W\dr^+(\dtilde)$ et $w_n \in W_e(\dtilde)$ et donc $\phi^{-n}(y-w_n) \in W\dr^+(\dtilde)$. On a alors :
\[ \phi^{-(n+1)}(y-w_{n+1}) - \phi^{-(n+1)}(\phi(y) - \phi(w_n)) = \phi^{-(n+1)}(x) - (w_{n+1} - w_n) \in W\dr^+(\dtilde), \]
ce qui fait (comme $\phi^{-(n+1)}(x)\in W\dr^+(\dtilde)$) que $w_{n+1} - w_n \in W\dr^+(\dtilde) \cap W_e(\dtilde) = \dtilde^{\phi=1}$ par la proposition \ref{xzdtphi}. Quitte à modifier $w_{n+1}$ par un élément de $\dtilde^{\phi=1}$, on peut donc supposer que $w_n$ ne dépend pas de $n$, ce qui fait qu'il existe $w \in (\dtilde[1/t])^{\phi=1}$ tel que $\phi^{-n}(y-w) \in W\dr^+(\dtilde)$ pour tout $n \gg 0$. Si l'on écrit $y-w = \sum_{i=1}^d b_i \otimes f_i$ avec $b_i \in \btrig{}[1/t]$ comme dans la démonstration de la proposition \ref{xzdtphi}, on trouve alors que $\phi^{-n}(b_i) \in \bdr^+$ pour tout $n \gg 0$ ce qui fait que $y-w \in \dtilde$ et donc que $x=(1-\phi)y \in (1-\phi)\dtilde$ et donc que $\alpha$ est injective. En revanche, le fait que $\alpha$ est surjective est moins évident et suit des calculs sous-jacents à la proposition \ref{prisoc1}.
\end{rema}

\begin{prop}\label{recslope}
Soit $\dtilde$ un $(\phi,G_K)$-module sur $\btrig{}$, $s$ sa plus grande pente, $k \geq 1$ et $b \in \ZZ$.
\begin{enumerate}
\item Il existe deux $(\phi,G_K)$-modules $\dtilde^1$ et $\dtilde^2$ sur $\btrig{}$ qui sont isoclines de pente $s$ et tels que $\dtilde = (\dtilde_{<s} \oplus \dtilde^1) / \dtilde^2$;
\item le $\Qp$-espace vectoriel $\dtilde^{\phi^k=p^b}$ est une presque $\Cp$-représentation.
\end{enumerate}
\end{prop}

\begin{proof}
La démonstration se fait par récurrence sur le nombre $\ell(\dtilde)$ de pentes de $\dtilde$. Si $\ell(\dtilde)=1$, alors $\dtilde$ est isocline de pente $s=a/h$, et pour montrer le (1) il suffit de prendre $\dtilde^1 = \dtilde$ et $\dtilde^2 = 0$ et le (2) suit de la proposition \ref{modisoc} qui dit que $\dtilde = \btrig{} \otimes_{\Qph} V_{a,h}$ et de la proposition \ref{prisoc0} qui dit que $(\btrig{} \otimes_{\Qph} V_{a,h})^{\phi^k=p^b}$ est une presque $\Cp$-représentation.

Sinon, soit $s=a/h$ la plus grande pente de $\dtilde$ ce qui fait que l'on a une suite exacte $0 \to \dtilde_{<s} \to \dtilde \to \dtilde_s \to 0$ où $\ell(\dtilde_{<s}) = \ell(\dtilde) -1$. En prenant les invariants par $\phi^h-p^a$, on trouve :
\[ 0 \to \dtilde_{<s}^{\phi^h=p^a} \to \dtilde^{\phi^h=p^a} \to \dtilde_s^{\phi^h=p^a} \to \dtilde_{<s} / (\phi^h - p^a). \]

Par le (2) de la proposition \ref{hiphi}, on a $\dtilde_{<s} / (\phi^h - p^a) = 0$; par ailleurs, $\dtilde_{<s}^{\phi^h=p^a}$ est une presque $\Cp$-représentation par hypothèse de récurrence et enfin par la proposition \ref{modisoc}, on peut écrire $\dtilde_s = \btrig{} \otimes_{\Qph} V_{a,h}$ et alors  $\dtilde_s^{\phi^h=p^a} = V_{a,h}$ est de dimension finie sur $\Qp$. On peut alors appliquer le (2) de la proposition \ref{extpcp} qui nous fournit $X \subset \dtilde^{\phi^h=p^a}$ de dimension finie qui s'envoie surjectivement sur $V_{a,h}$. Quitte à remplacer $X$ par le $\Qph$-espace vectoriel qu'il engendre, on peut supposer que $X$ est un $\Qph$-espace vectoriel. Posons alors $\dtilde^1 = \btrig{} \otimes_{\Qph} X$; il est isocline de pente $s=a/h$ puisque $X \subset \dtilde^{\phi^h=p^a}$ et on a une application naturelle : $\dtilde_{<s} \oplus \dtilde^1 \to \dtilde$ qui est surjective. Son noyau $\dtilde^2$ s'injecte dans $\dtilde^1$ et ses pentes sont donc $\geq s$ (par le (a) de la proposition 4.5.14 de \cite{KK05}). Par ailleurs, le lemme 3.4.3 de \cite{KK05} appliqué à la suite $0 \to \dtilde^2 \to \dtilde_{<s} \oplus \dtilde^1 \to \dtilde \to 0$ implique que :
\[ \deg(\dtilde^2) = \deg(\dtilde_{<s}) + \deg(\dtilde^1) - \deg(\dtilde) = \deg(\dtilde^1) - \deg(\dtilde) = s \cdot \rg(\dtilde^2), \]
ce qui fait que $\dtilde^2$ est nécessairement isocline de pente $s$, ce qui montre le (1).

Pour montrer le (2), posons $r=b/k$ et considérons la suite exacte :
\[ 0 \to \dtilde_{\leq r} \to \dtilde \to \dtilde_{>r} \to 0. \]
En prenant les invariants par $\phi^k-p^b$, on trouve $0 \to \dtilde_{\leq r}^{\phi^k=p^b} \to \dtilde^{\phi^k=p^b} \to \dtilde_{>r}^{\phi^k=p^b}$ et par le (1) de la proposition \ref{hiphi}, on a $\dtilde_{>r}^{\phi^k=p^b} = 0$ ce qui fait que l'application naturelle $\dtilde_{\leq r}^{\phi^k=p^b} \to \dtilde^{\phi^k=p^b}$ est un isomorphisme. On peut donc supposer que les pentes de $\dtilde$ sont $\leq r$. Dans ce cas, le (1) nous donne une suite exacte :
\[ 0 \to \dtilde^2 \to \dtilde_{<s} \oplus \dtilde^1 \to \dtilde \to 0, \]
et en prenant les invariants par $\phi^k-p^b$, on trouve :
\[ 0 \to (\dtilde^2)^{\phi^k=p^b} \to \dtilde_{<s}^{\phi^k=p^b} \oplus (\dtilde^1)^{\phi^k=p^b} \to \dtilde^{\phi^k=p^b} \to \dtilde^2 / (\phi^k-p^b), \]
Comme $s\leq r$, le (2) de la proposition 1.1.4 nous dit que $\dtilde^2 / (\phi^k-p^b) = 0$. Par ailleurs, $(\dtilde^1)^{\phi^k=p^b}$ et $(\dtilde^2)^{\phi^k=p^b}$ et $\dtilde_{<s}^{\phi^k=p^b}$ sont trois presque $\Cp$-représentations par hypothèse de récurrence ce qui fait que $\dtilde^{\phi^k=p^b}$ est elle-même une presque $\Cp$-représentation.
\end{proof}

\begin{proof}[Démonstration du théorème \ref{xipcp}]
Le (1) fait l'objet des propositions \ref{xzdtphi} et \ref{xundtsurphi}. Montrons le (2). Le fait que $X^0(W)$ est une presque $\Cp$-représen\-tation suit de la proposition \ref{xzdtphi} qui dit que $X^0(W) = \dtilde^{\phi=1}$ et du (2) de la proposition \ref{recslope} appliquée à $b/k=0/1$. 

Le fait que $X^1(W)$ est une presque $\Cp$-représen\-tation se démontre par récurrence sur le nombre de pentes de $\dtilde$. Si $\dtilde$ est isocline, alors le résultat suit des propositions \ref{modisoc} et \ref{prisoc1}.   Sinon, le (1) de la proposition \ref{recslope} nous fournit une suite exacte $0 \to \dtilde^2 \to \dtilde_{<s} \oplus \dtilde^1 \to \dtilde \to 0$ où $\dtilde^1$ et $\dtilde^2$ et $\dtilde_{<s}$ ont chacun moins de pentes que $\dtilde$ ce qui permet de montrer par récurrence que $X^1(\dtilde)$ est une presque $\Cp$-représentation en utilisant la suite exacte \ref{longxi} :
\[ 0 \to X^0(\dtilde^2) \to X^0(\dtilde_{<s} \oplus \dtilde^1) \to X^0(\dtilde) \to X^1(\dtilde^2) \to X^1(\dtilde_{<s} \oplus \dtilde^1) \to X^1(\dtilde) \to 0. \]

Pour montrer le (3), observons que d'une part, l'application $X^0(\dtilde_{\leq 0}) \to X^0(\dtilde)$ est un isomorphisme et d'autre part si $s \leq 0$, alors $X^1(\dtilde_{< s})=0$ et on a donc une suite exacte :
\[ 0 \to X^0(\dtilde_{< s}) \to X^0(\dtilde) \to X^0(\dtilde_s) \to 0. \]
La dimension étant additive sur les suites exactes, on se ramène au cas où $\dtilde$ est isocline, qui résulte alors de la proposition \ref{modisoc} et de la proposition \ref{prisoc0}. Le (4) se démontre exactement de la même manière, en utilisant la proposition \ref{prisoc1}. Ceci termine la démonstration du théorème.
\end{proof}

\begin{rema}\label{xandslope}
Notons pour référence les propriétés suivantes de $X^0$ et $X^1$ : 
\begin{enumerate}
\item les applications $X^0(\dtilde_{\leq 0}) \to X^0(\dtilde)$ et $X^1(\dtilde) \to X^1(\dtilde_{>0})$ sont des isomorphismes;
\item on a $X^0(\dtilde)=0$ si et seulement si les pentes de $\dtilde$ sont toutes $>0$ et on a $X^1(\dtilde)=0$ si et seulement si les pentes de $\dtilde$ sont toutes $\leq 0$.
\end{enumerate}
\end{rema}

\section{Pleine fidélité pour les $\be$-représentations}\label{parffbe}

L'objet de ce  pararaphe est de démontrer le théorème suivant, qui répond par l'affirmative à une question de Fontaine (voir la remarque en bas de la page 375 de \cite{F03}).

\begin{theo}\label{ffbe}
Le foncteur d'oubli de la catégorie des $\be$-représentations de $G_K$ vers la catégorie des $\Qp$-espaces vectoriels topologiques avec action linéaire et continue de $G_K$ est pleinement fidèle.
\end{theo}

Avant de montrer ce théorème, définissons la topologie naturelle sur une $\be$-représen\-tation. Rappelons que si $W_e$ est une $\be$-représentation et que si l'on pose $W\dr = \bdr \otimes_{\be} W_e$, alors (cf.\ le début du \S 3.5 de \cite{F04}) la $\bdr$-représentation $W\dr$ admet un $\bdr^+$-réseau $W\dr^+$ stable par $G_K$. Si $j \geq 0$, alors $W_e^j = W_e \cap t^{-j} W\dr^+$ est, par les résultats du paragraphe précédent, une presque $\Cp$-représentation, et c'est notamment un objet de $\bangk$. Comme on a $W_e = \cup_{j \geq 0} W_e^j$, la $\be$-représentation $W_e$ est une limite inductive dénombrable d'espaces de Banach et donc un espace LF. Par ailleurs, si l'on choisit un $\bdr^+$-réseau différent de $W\dr$, alors comme deux tels réseaux sont commensurables, on trouve une structure d'espace LF homéomorphe à la première et la topologie de $W_e$ qui en résulte ne dépend donc pas du choix de $W\dr^+$. La topologie d'une $\be$-représentation est alors la topologie d'espace LF que l'on a définie ci-dessus. Dans le cas où $W_e = \be \otimes_{\Qp} V$, on retrouve la définition du début du \S 8.2 de \cite{F03}.

Rappelons que si $f : E \to F$ est une application linéaire continue d'espaces LF avec $E= \cup_{i \geq 1} E_i$ et $F = \cup_{j \geq 1} F_j$, alors par le théorème de Baire, pour tout $i$ il existe $j$ tel que $f(E_i) \subset F_j$ et que réciproquement, une collection compatible de telles applications définit une application continue d'espaces LF. En particulier, si $f : W_e \to X_e$ est un morphisme de $\be$-représentations, alors il est nécessairement continu et le foncteur d'oubli $\rep_{\be}(G_K) \to \rep_{LF}(G_K)$ est bien défini et il est manifestement fidèle ce qui fait que pour montrer le théorème, il suffit de montrer que toute application continue $G_K$-équivariante $W_e \to X_e$ est nécessairement $\be$-linéaire.

\begin{proof}[Démonstration du théorème \ref{ffbe}]
Si $f : W_e \to X_e$ est une application continue $G_K$-équivariante, alors pour tout $i \geq 0$, il existe $j$ tel que $f(W_e \cap t^{-i} W\dr^+) \subset X_e \cap t^{-j} X\dr^+$. Quitte à modifier $W\dr^+$ et $X\dr^+$, on peut supposer que $f(W_e \cap W\dr^+) \subset X_e \cap X\dr^+$ et que les pentes de $(W_e,W\dr^+)$ sont $\leq 0$. Cette dernière hypothèse implique que les pentes de $(W_e,t^{-i} W\dr^+)$ sont $\leq -i$ (rappelons que si $\dfont(W)$ est le $(\phi,\Gamma)$-module associé à une $B$-paire $W= (W_e,W\dr^+)$, alors le $(\phi,\Gamma)$-module associé à $W= (W_e, t^j W\dr^+)$ est $t^j \dfont(W)$) et en particulier que si $i \geq 0$, alors :
\[ \dim_{\prcp}(W_e \cap t^{-i} W\dr^+) - \dim_{\prcp} (W_e \cap W\dr^+) = i \cdot \rg(W) \]
et $\haut(W_e \cap t^{-i} W\dr^+)=\haut(W_e \cap W\dr^+)$ par le (3) du théorème \ref{xipcp} ce qui fait que l'on a une suite exacte :
\[ 0 \to W_e \cap W\dr^+ \to W_e \cap t^{-i} W\dr^+ \to t^{-i} W\dr^+ / W\dr^+ \to 0, \]
et que si l'on choisit $j$ tel que $f(W_e \cap t^{-i} W\dr^+) \subset X_e \cap t^{-j} X\dr^+$, alors $f$ induit une application continue et $G_K$-équivariante $f\dr^{ij} : t^{-i} W\dr^+ / W\dr^+ \to t^{-j} X\dr^+ / X\dr^+$ puisque l'on a un diagramme :
\[ \begin{CD} 0 @>>> W_e \cap W\dr^+ @>>> W_e \cap t^{-i} W\dr^+ @>>> t^{-i} W\dr^+ / W\dr^+ @>>> 0 \\
&& @V{f}VV @V{f}VV @VVV \\\
0 @>>> X_e \cap X\dr^+ @>>> X_e \cap t^{-j} X\dr^+ @>>> t^{-j} X\dr^+ / X\dr^+.  
\end{CD} \]
Le théorème A' de \cite{F03} implique que $f\dr^{ij}$ est $\bdr^+$-linéaire. On trouve que $f$ induit une application $\bdr^+$-linéaire $f\dr^i : t^{-i} W\dr^+ / W\dr^+ \to X\dr / X\dr^+$ puis en passant à la limite une application $\bdr^+$-linéaire : $f\dr : W\dr / W\dr^+ \to X\dr / X\dr^+$. Cette application induit elle-même une application $\bdr^+$-linéaire $W\dr^+ \to X\dr^+$ (en passant au duaux) et donc une application $\bdr$-linéaire $W\dr \to X\dr$; comme $W_e \subset W\dr$ et $X_e \subset X\dr$ on trouve bien que $f$ est $\be$-linéaire.
\end{proof}

\appendix

\section{Réseaux et filtrations}\label{app}

Dans le \S III de \cite{Ber5} (la fin de la démonstration du théorème III.2.3) ainsi que dans le \S 3.3 de \cite{Ber8}  (la fin de la démonstration de la proposition 3.3.10), on utilise implicitement le résultat qui suit. Comme sa version sans l'action de $\Gamma$ est en fait fausse en général, il est utile de donner un énoncé correct et sa démonstration. Rappelons que si $\Gamma$ est un sous-groupe ouvert de $\Zp^\times$, on note $\nabla = \log(\gamma) / \log_p(\chi(\gamma))$ pour $\gamma \in \Gamma$ suffisamment proche de $1$.

\begin{prop}\label{misgam}
Soit $D$ un $K$-espace vectoriel muni d'une action de $\Gamma$ et $M$ un $K[\![t]\!]$-réseau stable par $\Gamma$ de $K(\!(t)\!) \otimes_K D$. 

Si l'on pose $\fil^i D = t^i M \cap D$, alors $M = \fil^0 (K(\!(t)\!) \otimes_K D)$, si l'on suppose qu'il existe $P(X) \in K[X]$ vérifiant $P(\nabla)=0$ sur $D$ et $P(X) \wedge P(X+j) = 1$ pour tout $j \neq 0$.
\end{prop}

\begin{proof}
Il est immédiat que $\fil^0 (K(\!(t)\!) \otimes_K D) \subset M$ et nous allons donc vérifier que l'on a bien $M \subset \fil^0 (K(\!(t)\!) \otimes_K D)$. Soit $e_1,\dots,e_d$ une base de $D$ adaptée à la filtration avec $e_i \in \fil^{h_i} D$  et $m \in M$. On peut écrire $m = \sum_{i=1}^d \alpha_i(t) t^{-h_i} e_i$ où $\alpha_i(t) = t^{-n_i} \alpha_i^*(t) \in t^{-n_i} K[\![t]\!]^\times$. Si $n_i \leq 0$ pour tout $i$, alors $m \in \fil^0 (K(\!(t)\!) \otimes_K D)$ et il n'y a rien à montrer. Dans le cas contraire, soit $n = \max_i n_i$; en écrivant $t^n m = \sum_{n=n_i} \alpha_i^*(t) t^{-h_i} e_i +  \sum_{n > n_i} t^{n-n_i} \alpha_i^*(t) t^{-h_i} e_i$, on trouve que $\sum_{n=n_i} \alpha_i^*(t) t^{-h_i} e_i \in t M$ et en posant $x_i = \alpha_i^*(0)$, on en déduit l'existence d'une relation non triviale :
\[ t^{-h_1} x_1 + \cdots + t^{-h_d} x_d \in tM, \] 
avec $x_j \in \fil^{h_j} D \setminus \fil^{h_j+1} D$ pour tout $j$. Parmi toutes les relations de ce type, prenons-en une non triviale de longueur minimale. Si on choisit $k$ tel que $x_k \neq 0$, alors on a : 
\[ P(\nabla+h_k)(t^{-h_1} x_1 + \cdots + t^{-h_d} x_d) = t^{-h_1} P(\nabla-h_1+h_k) x_1 + \cdots + t^{-h_d} P(\nabla-h_d+h_k) x_d. \]
L'hypothèse selon laquelle $P(X) \wedge P(X+j) = 1$ pour tout $j \neq 0$ et le fait que la filtration de $D$ est stable par $\Gamma$ impliquent que si $j \neq 0$, alors $P(\nabla+j)$ est un iso\-morphisme de $ \fil^h D / \fil^{h+1} D$ et comme l'on obtient ci-dessus une relation de longueur inférieure, et donc nulle, la relation de départ était du type $t^{-h_k} x_k \in t M$ ce qui implique que $x_k \in \fil^{h_k+1} D$, absurde. On avait donc bien $n_i \leq 0$ pour tout $i$ et $m \in \fil^0 (K(\!(t)\!) \otimes_K D)$.
\end{proof}

\end{document}